\documentclass[11pt,a4paper]{amsart}
\usepackage[utf8]{inputenc}
\usepackage[english]{babel}
\usepackage{pst-grad} 
\usepackage{pst-plot} 
\usepackage{pstricks}
\usepackage{amsmath}
\usepackage{amssymb}
\usepackage{mathrsfs}
\usepackage{amsthm}
\usepackage[dvips]{graphicx}
\usepackage{epsfig}

\newcommand{\RR}{\mathbb R}

\newcommand{\TT}{\mathbb T}

\newcommand{\pat}{\partial_t}
\newcommand{\pax}{\partial_x}

\newcommand{\vertiii}[1]{{\left\vert\kern-0.25ex\left\vert\kern-0.25ex\left\vert #1 
    \right\vert\kern-0.25ex\right\vert\kern-0.25ex\right\vert}}

\usepackage[pdfauthor={...},pdftitle={...},bookmarks,colorlinks]{hyperref}

\newcounter{comentcount}
\newcounter{teocount}
\newtheorem{corollary}{Corollary}
\newtheorem{theorem}[teocount]{Theorem}

\newtheorem{remark}{Remark}

\title[]{Global solutions for a hyperbolic-parabolic system of chemotaxis}

\author[R. Granero-Belinch\'{o}n]{Rafael Granero-Belinch\'{o}n}
\email{granero@math.univ-lyon1.fr}
\address{Univ Lyon, Universit\'e Claude Bernard Lyon 1, CNRS UMR 5208, Institut Camille Jordan, 43 blvd. du 11 novembre 1918, F-69622 Villeurbanne cedex, France.}

\begin{document}

\begin{abstract}
We study a hyperbolic-parabolic model of chemotaxis in dimensions one and two. In particular, we prove the global existence of classical solutions in certain dissipation regimes. 
\end{abstract}

\maketitle

{\small
\tableofcontents}
\section{Introduction}
In this note we study the following system of partial differential equations
\begin{align}
\pat u&=-\Lambda^{\alpha} u+\nabla\cdot (uq),\text{ for }x\in\TT^d,\,t\geq0,\label{eq:6b}\\
\pat q&=\nabla f(u),\text{ for }x\in\TT^d,\,t\geq0,\label{eq:6b2}
\end{align}
where $u$ is a non-negative scalar function, $q$ is a vector in $\RR^d$, $\TT^d$ denotes the domain $[-\pi,\pi]^d$ with periodic boundary conditions, $d=1,2$ is the dimension, $f(u)=u^2/2$, $0<\alpha\leq 2$ and $(-\Delta)^{\alpha/2}=\Lambda^{\alpha}$ is the fractional Laplacian. 

This system was proposed by Othmers \& Stevens \cite{stevens1997aggregation} based on biological considerations as a model of tumor angiogenesis. In particular, in the previous system, $u$ is the density of vascular endothelial cells and $q=\nabla\log(v)$ where $v$ is the concentration of the signal protein known as vascular endothelial growth factor (VEGF) (see Bellomo, Li, \& Maini \cite{bellomo2008foundations} for more details on tumor modelling). Similar hyperbolic-dissipative systems arise also in the study of compressible viscous fluids or magnetohydrodynamics (see S. Kawashima \cite{kawashima1986large} and the references therein).

Equation \eqref{eq:6b} appears as a singular limit of the following Keller-Segel model of aggregation of the slime mold \emph{Dictyostelium discoideum} \cite{keller1970initiation} (see also Patlak \cite{patlak1953random})
\begin{equation}\label{eq:1}
\left\{\begin{aligned}
\pat u&=-\Lambda^\alpha u-\chi\nabla\cdot (u\nabla G(v)),\\
\pat v&=\nu\Delta v+(f(u)+\lambda)v,
\end{aligned}\right.
\end{equation}
when $G(v)=\log(v)$ and the diffusion of the chemical is negligible, \emph{i.e.} $\nu\rightarrow0$. 
 
Similar equations arising in different context are the Majda-Biello model of Rossby waves \cite{majda2003nonlinear} or the magnetohydrodynamic-Burgers system proposed by Fleischer \& Diamond \cite{fleischer2000burgers}.

Most of the results for \eqref{eq:6b} corresponds to the case where $d=1$. Then, when the diffusion is local \emph{i.e.} $\alpha=2$, \eqref{eq:6b} has been studied by many different research groups. In particular, Fan \& Zhao \cite{fan2012blow}, Li \& Zhao \cite{li2015initial}, Mei, Peng \& Wang \cite{mei2015asymptotic}, Li, Pan \& Zhao \cite{li2015quantitative}, Jun, Jixiong, Huijiang \& Changjiang \cite{jun2009global} Li \& Wang \cite{li2009nonlinear} and Zhang \& Zhu \cite{zhang2007global} studied the system \eqref{eq:6b} when $\alpha=2$ and $f(u)=u$ under different boundary conditions (see also the works by Jin, Li \& Wang \cite{jin2013asymptotic}, Li, Li \& Wang \cite{li2014stability},  Wang \& Hillen \cite{wang2008shock} and Wang, Xiang \& Yu \cite{wang2016asymptotic}). The case with general $f(u)$ was studied by Zhang, Tan \& Sun \cite{zhang2013global} and Li \& Wang \cite{li2010nonlinear}. 

Equation \eqref{eq:6b} in several dimensions has been studied by Li, Li \& Zhao \cite{li2011hyperbolic}, Hau \cite{hao2012global} and Li, Pan \& Zhao \cite{li2012global}. There, among other results, the global existence for small initial data in $H^s$, $s>2$ is proved.

To the best of our knowledge, the only result when the diffusion is nonlocal, \emph{i.e.} $0<\alpha<2$, is \cite{Ghyperparweak}. In that paper we obtained appropriate lower bounds for the fractional Fisher information and, among other results, we proved the global existence of weak solution for $f(u)=u^r/r$ and $1\leq r\leq 2$. 

In this note, we address the existence of classical solutions in the case $0<\alpha\leq 2$. This is a challenging issue due to the hyperbolic character of the equation for $q$. In particular, $u$ verifies a transport equation where the velocity $q$ is one derivative more singular than $u$ (so $\nabla \cdot (u q)$ is two derivatives less regular than $u$).

\section{Statement of the results}\label{sec2}
For the sake of clarity, let us state some notation: we define the mean as
$$
\langle g\rangle =\frac{1}{(2\pi)^d}\int_{\TT^d}g(x)dx.
$$
Also, from this point onwards, we write $H^s$ for the $L^2-$based Sobolev space of order $s$ endowed with the norm
$$
\|u\|_{H^s}^2=\|u\|_{L^2}^2+\|u\|_{\dot{H}^s}^2,\quad \|u\|_{\dot{H}^s}=\|\Lambda^s u\|_{L^2}.
$$
For $\beta\geq0$, we consider the following energies $E_\beta$ and dissipations $D_\beta$,
\begin{equation}\label{energy}E_\beta(t)=\|u\|_{\dot{H}^\beta}^2+\|q\|_{\dot{H}^\beta}^2,\quad D_\beta(t)=\|u\|_{\dot{H}^{\beta+\alpha/2}}^2.
\end{equation}
Recall that the lower order norms verify the following energy balance \cite{Ghyperparweak}
\begin{equation}\label{lowerorder}
\frac{1}{2}\left(\|u(t)\|_{L^2}^2+\|q(t)\|_{L^2}^2\right)+\int_0^t\|u(s)\|_{\dot{H}^{\alpha/2}}^2ds=\frac{1}{2}\left(\|u_0\|_{L^2}^2+\|q_0\|_{L^2}^2\right).
\end{equation}

\subsection{On the scaling invariance}\label{sec2e}
Notice that the equations \eqref{eq:6b}-\eqref{eq:6b2} verify the following scaling symmetry: for every $\lambda>0$
$$
u_\lambda(x,t)=\lambda^{\alpha-1} u\left(\lambda x,\lambda^\alpha t\right),\quad q_\lambda(x,t)=\lambda^{\alpha-1} q\left(\lambda x,\lambda^\alpha t\right).
$$
This scaling serves as a \emph{zoom in} towards the small scales. We also know that 
$$
\|u(t)\|_{L^2}^2+\|q(t)\|_{L^2}^2
$$
is the strongest (known) quantity verifying a global-in-time bound. Then, in the one dimensional case, the $L^2$ norms of $u$ and $q$ are invariant under the scaling of the equations when $\alpha=1.5$. That makes $\alpha=1.5$ the critical exponent for the global estimates known. Equivalently, if we define the rescaled (according to the scaling of the strongest conserved quantity $\|u(t)\|_{L^2}^2+\|q(t)\|_{L^2}^2$) functions  
$$
u_\gamma(x,t)=\gamma^{0.5} u\left(\gamma x,\gamma^\alpha t\right),\quad q_\gamma(x,t)=\gamma^{0.5} q\left(\gamma x,\gamma^\alpha t\right).
$$
we have that $u_\gamma$ and $q_\gamma$ solve
\begin{align*}
\pat u_\gamma&=-\Lambda^{\alpha} u_\gamma+\gamma^{\alpha-1.5}\pax (u_\gamma q_\gamma),\\
\pat q_\gamma&=\gamma^{\alpha-1.5}u_\gamma \pax u_\gamma.
\end{align*}

Larger values of $\alpha$ form the subcritical regime where the diffusion dominates the drift in small scales. Smaller values of $\alpha$ form the supercritical regime where the drift might be dominant at small scales.

Similarly, the two dimensional case has critical exponent $\alpha=2$.

\begin{remark}
Notice that the equations \eqref{eq:6b}-\eqref{eq:6b2} where $f(u)=u$ have a different scaling symmetry but the same critical exponent $\alpha=1.5$. In this case, the scaling symmetry is given by
$$
u_\lambda(x,t)=\lambda^{2\alpha-2} u\left(\lambda x,\lambda^\alpha t\right),\quad
q_\lambda(x,t)=\lambda^{\alpha-1}q\left(\lambda x,\lambda^\alpha t\right),
$$
while the conserved quantity is $\|u(t)\|_{L^1}+ \|q(t)\|_{L^2}^2/2.$
Thus, if we define the rescaled (according to the scaling of the conserved quantity) functions  
$$
u_\gamma(x,t)=\gamma u\left(\gamma x,\gamma^\alpha t\right),\quad q_\gamma(x,t)=\gamma^{0.5} q\left(\gamma x,\gamma^\alpha t\right).
$$
we have that $u_\gamma$ and $q_\gamma$ solve
\begin{align*}
\pat u_\gamma&=-\Lambda^{\alpha} u_\gamma+\gamma^{\alpha-1.5}\pax (u_\gamma q_\gamma),\\
\pat q_\gamma&=\gamma^{\alpha-1.5} \pax u_\gamma.
\end{align*}
A global existence result when $\alpha$ is the range $1.5\leq \alpha<2$ for the problem where $f(u)=u$ is left for future research.
\end{remark}

\subsection{Results in the one-dimensional case $d=1$}\label{sec2b}
One of our main results is
\begin{theorem}\label{globalstrong} Fix $T$ an arbitrary parameter and let $(u_0,q_0)\in H^2(\TT)\times H^2(\TT)$ be the initial data such that $0\leq u_0$ and $\langle q_0\rangle=0$. Assume that $\alpha\geq1.5$. Then there exist a unique global solution $(u(t),q(t))$ to problem \eqref{eq:6b} verifying
$$
u\in L^\infty(0,T;H^2(\TT))\cap L^2(0,T;H^{2+\alpha/2}(\TT)),
q\in L^\infty(0,T;H^{2}(\TT)).
$$
Furthermore, the solution is uniformly bounded in
$$
(u,q)\in C([0,\infty),H^1(\TT))\times C([0,\infty),H^1(\TT)).
$$
\end{theorem} 

In the case where the strength of the diffusion, $\alpha$, is even weaker, we have the following global existence result for small data:

\begin{theorem}\label{globalsmall} Fix $T$ an arbitrary parameter and let $(u_0,q_0)\in H^2(\TT)\times H^2(\TT)$ be the initial data such that $0\leq u_0$ and $\langle q_0\rangle=0$. There exists $\mathcal{C}_\alpha$ such that if $1.5>\alpha>1$ and 
$$
\|u_0\|_{\dot{H}^{\alpha/2}}^2+\|q_0\|_{\dot{H}^{\alpha/2}}^2\leq \mathcal{C}_\alpha
$$ 
then there exist a unique global solution $(u(t),q(t))$ to problem \eqref{eq:6b} verifying
$$
u\in L^\infty(0,T;H^2(\TT))\cap L^2(0,T;H^{2+\alpha/2}(\TT)),
q\in L^\infty(0,T;H^{2}(\TT)).
$$
Furthermore, the solution verifies
$$
\|u(t)\|_{\dot{H}^{\alpha/2}}^2+\|q(t)\|_{\dot{H}^{\alpha/2}}^2\leq \|u_0\|_{\dot{H}^{\alpha/2}}^2+\|q_0\|_{\dot{H}^{\alpha/2}}^2.
$$
\end{theorem} 

\begin{corollary}\label{globalsmall2} Fix $T$ an arbitrary parameter and let $(u_0,q_0)\in H^2(\TT)\times H^2(\TT)$ be the initial data such that $0\leq u_0$ and $\langle q_0\rangle=0$. Assume that $1\geq\alpha\geq0.5$ and 
$$
\|u_0\|_{\dot{H}^1}^2+\|q_0\|_{\dot{H}^1}^2< \frac{4}{9C^2_S}
$$ 
where $C_S$ is defined in \eqref{Sobolev}. Then there exist a unique global solution $(u(t),q(t))$ to problem \eqref{eq:6b} verifying
$$
u\in L^\infty(0,T;H^2(\TT))\cap L^2(0,T;H^{2+\alpha/2}(\TT)),
q\in L^\infty(0,T;H^{2}(\TT)).
$$
Furthermore, the solution verifies
$$
\|u(t)\|_{H^1}^2+\|q(t)\|_{H^1}^2\leq \|u_0\|_{H^1}^2+\|q_0\|_{H^1}^2.
$$
\end{corollary} 

\subsection{Results in the two-dimensional case $d=2$}\label{sec2c}
In two dimensions the global existence read
\begin{theorem}\label{globalstrongd2} Fix $T$ an arbitrary parameter and let $(u_0,q_0)\in H^2(\TT^2)\times H^2(\TT^2)$ be the initial data such that $0\leq u_0$, $\langle q_0\rangle=0$ and $\text{curl}\, q_0=0$. Assume that $\alpha=2$. Then there exist a unique global solution $(u(t),q(t))$ to problem \eqref{eq:6b} verifying
$$
u\in L^\infty(0,T;H^2(\TT^2))\cap L^2(0,T;H^{3}(\TT^2)),
q\in L^\infty(0,T;H^{2}(\TT^2)).
$$
Furthermore, the solution is uniformly bounded in
$$
(u,q)\in C([0,\infty),H^1(\TT^2))\times C([0,\infty),H^1(\TT^2)).
$$
\end{theorem} 

\begin{corollary}\label{globalsmall2d2} Fix $T$ an arbitrary parameter and let $(u_0,q_0)\in H^2(\TT^2)\times H^2(\TT^2)$ be the initial data such that $0\leq u_0$, $\langle q_0\rangle=0$ and $\text{curl}\, q_0=0$. Assume that $2>\alpha\geq1$ and 
$$
\|u_0\|_{\dot{H}^1}^2+\|q_0\|_{\dot{H}^1}^2< \mathcal{C}
$$ 
where $\mathcal{C}$ is a universal constant. Then there exist a unique global solution $(u(t),q(t))$ to problem \eqref{eq:6b} verifying
$$
u\in L^\infty(0,T;H^2(\TT^2))\cap L^2(0,T;H^{2+\alpha/2}(\TT^2)),
q\in L^\infty(0,T;H^{2}(\TT^2)).
$$
Furthermore, the solution verifies
$$
\|u(t)\|_{H^1}^2+\|\nabla\cdot q(t)\|_{L^2}^2\leq \|u_0\|_{H^1}^2+\|\nabla\cdot q_0\|_{L^2}^2.
$$
\end{corollary} 

\begin{remark}
In the case where the domain is the one-dimensional torus, $\TT$, local existence of solution for \eqref{eq:6b}-\eqref{eq:6b2} was proved in \cite{Ghyperparweak} for a more general class of kinetic function $f(u)$. The local existence of solution for \eqref{eq:6b}-\eqref{eq:6b2} the domain is the two-dimensional torus $\TT^d$ with $d=2$ follows from the local existence result in \cite{Ghyperparweak} with minor modifications. Consequently, we will focus on obtaining global-in-time \emph{a priori} estimates.
\end{remark}

\subsection{Discussion}\label{sec2d}
Due to the hyperbolic character of the equation for $q$, prior available global existence results of classical solution for equation \eqref{eq:6b} impose several assumptions. Namely,
\begin{itemize}
\item either $d=1$ and $\alpha=2$ \cite{zhang2013global, li2010nonlinear}, 
\item or $d=2,3$, $\alpha=2$ and the initial data verifies some smallness condition on \emph{strong} Sobolev spaces $H^s$, $s\geq2$ \cite{xie2013global, zhang2015global}.
\end{itemize}

Our results removed some of the previous conditions. On the one hand, we prove global existence for arbitrary data in the cases $d=1$ and $\alpha\geq1.5$ and $d=2$ and $\alpha=2$. On the other hand, in the cases where we have to impose size restrictions on the initial data, the Sobolev spaces are bigger than $H^2$ (thus, the norm is weaker). Finally, let us emphasize that our results can be adapted to the case where the spatial domain is $\RR^d$.

A question that remains open is the trend to equilibrium. From \eqref{lowerorder} is clear that the solution $(u(t),q(t))$ tends to the homogeneous state, namely $(\langle u_0 \rangle,0)$. However, the rate of this convergence is not clear.

\section{Proof of Theorem \ref{globalstrong}}
\textbf{Step 1; $H^1$ estimate:} Testing the first equation in \eqref{eq:6b} against $\Lambda^2 u$, integrating by parts and using the equation for $q$, we obtain
\begin{align*}
\frac{1}{2}\frac{d}{dt}\|u\|_{\dot{H}^1}^2+\|\pax u\|_{\dot{H}^{\alpha/2}}^2&=-\int_\TT \pax(uq)\pax^2udx\\
&=\frac{1}{2}\int_\TT \pax q (\pax u)^2dx-\int_\TT \pax q u\pax^2udx\\
&=\frac{1}{2}\int_\TT \pax q (\pax u)^2dx-\int_\TT \pax q (\pat \pax q -(\pax u)^2)dx,
\end{align*}
so
$$
\frac{1}{2}\frac{d}{dt}(\|u\|_{\dot{H}^1}^2+\|q\|_{\dot{H}^1}^2)+\|u\|_{\dot{H}^{1+\alpha/2}}^2=\frac{3}{2}\int_\TT \pax q (\pax u)^2dx.
$$
Denoting
$$
I=\frac{3}{2}\int_\TT \pax q (\pax u)^2dx,
$$
and using Sobolev embedding and interpolation, we have that
\begin{equation}
I\leq \frac{3}{2}\|q\|_{\dot{H}^1}\|\pax u\|_{L^4}^2\leq \frac{3}{2} C_S\|q\|_{\dot{H}^1}\|\pax u\|_{\dot{H}^{0.25}}^2\label{eqI},
\end{equation}
where $C_S$ is the constant appearing in the embedding
\begin{equation}\label{Sobolev}
\|g\|_{L^4}\leq C_S\|g\|_{\dot{H}^{0.25}}.
\end{equation}
Using the interpolation 
$$
H^{1+\alpha/2} \subset H^{1.25}\subset H^{\alpha/2},
$$
and Poincar\'e inequality (if $\alpha>1.5$) we conclude
$$
I\leq c\|q\|_{\dot{H}^1}\|\Lambda^{\alpha/2}u\|_{L^2}\|u\|_{\dot{H}^{1+\alpha/2}},
$$

Using \eqref{energy}, we have that
$$
\frac{d}{dt}E_1+D_1\leq c \|u\|_{\dot{H}^{\alpha/2}}^2 E_1.
$$
Using Gronwall's inequality and the estimate \eqref{lowerorder}, we have that
$$
\sup_{0\leq t<\infty}E_1(t)\leq C(\|u_0\|_{H^1},\|q_0\|_{H^1}),
$$
$$
\int_0 ^T D_1(s)ds\leq C(\|u_0\|_{H^1},\|q_0\|_{H^1}, T),\,\forall 0<T<\infty.
$$

\textbf{Step 2; $H^2$ estimate:} Now we prove that the solutions satisfying the previous bounds for $E_1$ and $D_1$ also satisfy the corresponding estimate in $H^2$. We test the equation for $u$ against $\Lambda^4 u$. We have that
\begin{align*}
\frac{1}{2}\frac{d}{dt}\|u\|_{\dot{H}^2}^2+\|u\|_{\dot{H}^{2+\alpha/2}}^2&=-\int_\TT \pax^2(uq)\pax^3udx\\
&=\int_\TT \pax q \frac{5(\pax^2 u)^2}{2}dx-\int_\TT \pax^2 q (\pat \pax^2 q -5\pax u\pax^2 u)dx,
\end{align*}
so
$$
\frac{1}{2}\frac{d}{dt}(\|u\|_{\dot{H}^2}^2+\|q\|_{\dot{H}^2}^2)+\|u\|_{\dot{H}^{2+\alpha/2}}^2=\frac{5}{2}\int_\TT \pax q (\pax^2 u)^2dx+5\int_\TT \pax^2 q \pax^2 u \pax udx.
$$
We define
$$
J_1=\frac{5}{2}\int_\TT \pax q (\pax^2 u)^2dx,\;
J_2=5\int_\TT \pax^2 q \pax^2 u \pax udx.
$$
Then, we have that
$$
J_1\leq c\|\pax q\|_{L^\infty}\|\pax^2 u\|_{L^2}^2\leq c\|\pax^2 q\|_{L^2}^{0.5}\|u\|_{\dot{H}^{1+\alpha/2}}^{\alpha}\|u\|_{\dot{H}^{2+\alpha/2}}^{2-\alpha},
$$
so, using Young's inequality,
$$
J_1\leq c\|\pax^2 q\|_{L^2}^{\frac{1}{\alpha}}\|u\|_{\dot{H}^{1+\alpha/2}}^{2}+\frac{1}{4}\|u\|_{\dot{H}^{2+\alpha/2}}^{2}.
$$
Similarly, using Poincar\'e inequality and $\alpha\geq0.5$,
$$
J_2\leq c\|\pax u\|_{L^4}\|\pax^2 u\|_{L^4}\|\pax^2 q\|_{L^2}\leq c\|u\|_{\dot{H}^{1+\alpha/2}}\|u\|_{\dot{H}^{2+\alpha/2}}\|\pax^2 q\|_{L^2},
$$
and
$$
J_2\leq c\|u\|_{\dot{H}^{1+\alpha/2}}^{2}\|\pax^2 q\|_{L^2}^2+\frac{1}{4}\|u\|_{\dot{H}^{2+\alpha/2}}^{2}.
$$
Finally,
$$
\frac{d}{dt}E_{2}(t)+D_{2}(t)\leq c\|u\|_{\dot{H}^{1+\alpha/2}}^{2}(E_{2}(t)+1)
$$
and we conclude using Gronwall's inequality.

\section{Proof of Theorem \ref{globalsmall}}
\textbf{Step 1; $H^{\alpha/2}$ estimate:} Testing the first equation in \eqref{eq:6b} against $\Lambda^{\alpha} u$, we obtain
\begin{align*}
\frac{1}{2}\frac{d}{dt}\|u\|_{\dot{H}^{\alpha/2}}^2+\|u\|_{\dot{H}^{\alpha}}^2&=\int_\TT \pax(uq)\Lambda^{\alpha} udx\\
&=-\int_\TT \Lambda^\alpha(uq)\pax udx\\
&=-\int_\TT \left(\Lambda^\alpha(uq)-u\Lambda^\alpha q\right)\pax udx-\int_\TT \Lambda^\alpha q u\pax udx,
\end{align*}
so
$$
\frac{1}{2}\frac{d}{dt}\left(\|u\|_{\dot{H}^{\alpha/2}}^2+\|q\|_{\dot{H}^{\alpha/2}}^2\right)+\|u\|_{\dot{H}^{\alpha}}^2\leq -\int_\TT [\Lambda^\alpha,u]q\pax udx.
$$
We define
$$
K=-\int_\TT [\Lambda^\alpha,u]q\pax udx.
$$
Using the classical Kenig-Ponce-Vega commutator estimate \cite{KenigPonceVega} and Sobolev embedding, we have that
\begin{align}
\|[\Lambda^\alpha,u]q\|_{L^2}&\leq c\left(\|\pax u\|_{L^{2+\epsilon}}\|\Lambda^{\alpha-1}q\|_{L^{\frac{4+2\epsilon}{\epsilon}}}+\|\Lambda^\alpha u\|_{L^{2}}\|q\|_{L^{\infty}}\right)\nonumber\\
&\leq c\left(\|u\|_{\dot{H}^{1+\frac{\epsilon}{4+2\epsilon}}}\|\Lambda^{\alpha-1}q\|_{\dot{H}^{\frac{1}{2}-\frac{\epsilon}{4+2\epsilon}}}+\|u\|_{\dot{H}^{\alpha}}\|q\|_{\dot{H}^{\frac{\alpha}{2}}}\right)\label{KPV}.
\end{align}
Thus, taking $\epsilon$ such that
$$
1+\frac{\epsilon}{4+2\epsilon}=\alpha,\;\text{ \emph{i.e.} }\epsilon=\frac{4\alpha-4}{3-2\alpha}
$$
Equation \eqref{KPV} reads
\begin{equation}\label{KPV2}
\|[\Lambda^\alpha,u]q\|_{L^2}\leq c\left(\|u\|_{\dot{H}^{\alpha}}\|q\|_{\dot{H}^{\frac{1}{2}}}+\|u\|_{\dot{H}^{\alpha}}\|q\|_{\dot{H}^{\frac{\alpha}{2}}}\right)
\end{equation}
Using \eqref{KPV2} and Poincar\'e inequality, we have that
\begin{align*}
K&\leq c\|u\|_{\dot{H}^{\alpha}}^2\|q\|_{\dot{H}^{\frac{\alpha}{2}}}
\end{align*}
Then, we have that
$$
\frac{d}{dt}E_{\frac{\alpha}{2}}+D_{\frac{\alpha}{2}}\leq c\sqrt{E_{\frac{\alpha}{2}}}D_{\frac{\alpha}{2}}.
$$
Thus, due to the smallness restriction on the initial data, we obtain
$$
E_{\frac{\alpha}{2}}(t)+\delta\int_0^t D_{\frac{\alpha}{2}}(s)ds\leq E_{\frac{\alpha}{2}}(0)
$$
for $0<\delta$ small enough. 

\textbf{Step 2; $H^{1}$ estimate:} Our starting point is \eqref{eqI}. Then we use the interpolation
$$
\|g\|^{2}_{\dot{H}^{0.25}}\leq c\|g\|_{L^2}\|g\|_{\dot{H}^{0.5}}, 
$$
to obtain
$$
I\leq c\|q\|_{\dot{H}^1}\|u\|_{\dot{H}^1}\|u\|_{\dot{H}^{1.5}}\leq cE_1D_{\frac{\alpha}{2}}+\frac{D_1}{2}.
$$
Collecting all the estimates, we have that
$$
\frac{d}{dt}E_{1}+D_{1}\leq cE_1D_{\frac{\alpha}{2}},
$$
and we conclude using Gronwall's inequality. The $H^2$ estimates follows as in the proof of Theorem \ref{globalstrong}.

\section{Proof of Corollary \ref{globalsmall2}}
Using $\alpha\geq0.5$ and the estimate \eqref{eqI}, we have that
$$
I\leq \frac{3}{2}C_S\|q\|_{\dot{H}^1}\|u\|_{\dot{H}^{1+\alpha/2}}^2\leq \frac{3}{2}C_S \sqrt{E_1} D_1.
$$
Thus,
$$
\frac{1}{2}\frac{d}{dt}E_1+D_1\leq \frac{3}{2}C_S \sqrt{E_1} D_1.
$$
Thus, due to the smallness restriction on the initial data, we obtain
$$
E_1(t)+\delta\int_0^t D_1(s)ds\leq E_1(0)
$$
for $0<\delta$ small enough. Equipped with this estimates, we can repeat the argument as in Step 2 in Theorem \ref{globalstrong}.

\section{Proof of Theorem \ref{globalstrongd2}}
Recall that the condition 
$$
\text{curl}\,q_0=0
$$
propagates in time, \emph{i.e.}
$$
\text{curl}\,q(t)=\text{curl}\,q_0+\frac{1}{2}\int_0^t \text{curl}\nabla u^2ds=0.
$$
Using Plancherel Theorem, we have that
\begin{align*}
\|\nabla q\|_{L^2}^2&=C\sum_{\xi\in \mathbb{Z}^2}|\xi|^2|\hat{q}(\xi)|^2\\
&=C\sum_{\xi\in \mathbb{Z}^2}(\xi_1^2+\xi_2^2)(\hat{q}_1^2+\hat{q}_2^2).
\end{align*}
Due to the irrotationality
$$
\xi^\perp\cdot \hat{q}=0.
$$
Then, we compute
\begin{align*}
\|\nabla \cdot q\|_{L^2}^2&=C\sum_{\xi\in \mathbb{Z}^2}|\xi\cdot\hat{q}(\xi)|^2\\
&=C\sum_{\xi\in \mathbb{Z}^2}(\xi_1\hat{q}_1(\xi)+\xi_2\hat{q}_2(\xi))^2\\
&=C\sum_{\xi\in \mathbb{Z}^2}(\xi_1\hat{q}_1(\xi))^2+(\xi_2\hat{q}_2(\xi))^2+2\xi_1\hat{q}_1(\xi)\xi_2\hat{q}_2(\xi)\\
&=C\sum_{\xi\in \mathbb{Z}^2}(\xi_1\hat{q}_1(\xi))^2+(\xi_2\hat{q}_2(\xi))^2+(\xi_2\hat{q}_1(\xi))^2+(\xi_1\hat{q}_2(\xi))^2.
\end{align*}
So, the vector field $q$ satisfies
$$
\|\nabla q\|_{L^2}\leq \|\nabla\cdot q\|_{L^2}.
$$
As a consequence of $\langle \pat q_i \rangle=0$ and $\langle q_0\rangle=0$, every coordinate of $q$ satisfy
$$
\langle q_i(t)\rangle=0,
$$
and the Poincar\'e-type inequality
\begin{equation}\label{eqpoinc}
\|q\|_{L^2}\leq c\|\nabla\cdot q\|_{L^2}.
\end{equation}

Notice that in two dimensions we also have the energy balance \eqref{lowerorder}. We test equation \eqref{eq:6b} against $\Lambda^2u$ and use the equation for $q$. We obtain
$$
\frac{1}{2}\frac{d}{dt}\left(\|u\|_{\dot{H}^1}^2+\|\nabla\cdot{q}\|_{L^2}^2\right)=-\|u\|^2_{\dot{H}^{2}}-\int_{\TT^2}\nabla u\cdot q\Delta udx+\int_{\TT^2}|\nabla u|^2\nabla\cdot qdx.
$$
Using H\"{o}lder inequality, Sobolev embedding and interpolation, we have that
\begin{align*}
\frac{d}{dt}\left(\|u\|_{\dot{H}^1}^2+\|\nabla\cdot{q}\|_{L^2}^2\right)+2\|u\|^2_{\dot{H}^{2}}&\leq c\left(\|u\|_{\dot{H}^{1.5}}\|q\|_{L^4}\|u\|_{\dot{H}^2}+\|u\|_{\dot{H}^{1.5}}^2\|\nabla\cdot q\|_{L^2}\right)\\
&\leq c\|u\|_{\dot{H}^{1}}^{0.5}\|q\|_{L^2}^{0.5}\|q\|_{H^1}^{0.5}\|u\|_{\dot{H}^2}^{1.5}\\
&+c\|u\|_{\dot{H}^{1}}\|u\|_{\dot{H}^{2}}\|\nabla\cdot q\|_{L^2}.
\end{align*}
Using the H\"{o}dge decomposition estimate together with the irrotationality of $q$ and \eqref{eqpoinc}, we have that 
\begin{equation}\label{eqaux}
\|q\|_{H^1}\leq c\left(\|q\|_{L^2}+\|\nabla\cdot q\|_{L^2}\right)\leq c\|\nabla\cdot q\|_{L^2}.
\end{equation}
Due to \eqref{lowerorder}, we obtain that
$$
\|q\|_{L^\infty(0,\infty,L^2)}^2+\|u\|_{L^2(0,\infty,\dot{H}^1)}^2\leq C(\|u_0\|_{L^2},\|q_0\|_{L^2})
$$
so,
\begin{align*}
\frac{d}{dt}\left(\|u\|_{\dot{H}^1}^2+\|\nabla\cdot{q}\|_{L^2}^2\right)+\|u\|^2_{\dot{H}^{2}}&\leq c\|u\|_{\dot{H}^{1}}^{2}\|\nabla\cdot q\|_{L^2}^{2}.
\end{align*}
Using Gronwall's inequality and the integrability of $\|u\|_{\dot{H}^1}^2$ (see \eqref{lowerorder}), we obtain 
$$
E_1\leq C(\|u_0\|_{H^1},\|q_0\|_{H^1}),
$$
$$
\int_0 ^T D_1(s)ds\leq C(\|u_0\|_{H^1},\|q_0\|_{H^1}, T),\,\forall 0<T<\infty.
$$
To obtain the $H^2$ estimates, we test against $\Lambda^4 u$. Then, using the previous $H^1$ uniform bound and
$$
\|q\|_{L^\infty}^2\leq c\|q\|_{L^2}\|q\|_{H^2}\leq c\|q\|_{L^2}\|\Delta q\|_{L^2},
$$
we have that
\begin{multline}\label{eqa}
\frac{1}{2}\frac{d}{dt}\|\Delta u\|_{L^2}^2+\|u\|_{\dot{H}^3}^2=-\int_{\TT^d} \nabla \Delta u \nabla(\nabla u\cdot q)dx\\
-\int_{\TT^d} u\nabla \Delta u \cdot \nabla(\nabla\cdot q)dx
-\int_{\TT^d} \nabla u\cdot\nabla \Delta u \nabla\cdot qdx.
\end{multline}
Due to the irrotationality of $q$ and the identity
$$
\nabla\nabla\cdot q-\Delta q=\text{curl}\,(\text{curl}\,q),
$$ 
we have
$$
\pat \Delta q=\pat \nabla(\nabla\cdot q)=\nabla|\nabla u|^2+u\nabla\Delta u+\nabla u\Delta u.
$$
Applying Sobolev embedding and interpolation, we obtain that \eqref{eqa} can be estimated as
$$
\frac{d}{dt}\left(\|\Delta u\|_{L^2}^2+\|\Delta q\|_{L^2}^2\right)+\|u\|_{\dot{H}^3}^2\leq c\|\Delta q\|_{L^2}^2,
$$
so,
$$
\frac{d}{dt}E_2+D_2\leq cE_2,
$$
and we conclude using Gronwall's inequality.

\section{Proof of Corollary \ref{globalsmall2d2}}
We test the equation \eqref{eq:6b} against $\Lambda^2 u$. We obtain that
$$
\frac{1}{2}\frac{d}{dt}\left(\|u\|_{\dot{H}^1}^2+\|\nabla\cdot{q}\|_{L^2}^2\right)=-\|u\|^2_{\dot{H}^{1+\frac{\alpha}{2}}}+\int_{\TT^2}\nabla(\nabla u\cdot q)\nabla udx+\int_{\TT^2}|\nabla u|^2\nabla\cdot qdx.
$$
After a short computation, using H\"{o}lder estimates, Sobolev embedding and interpolation, we obtain that
$$
\frac{1}{2}\frac{d}{dt}\left(\|u\|_{\dot{H}^1}^2+\|\nabla\cdot{q}\|_{L^2}^2\right)+\|u\|^2_{\dot{H}^{1+\frac{\alpha}{2}}}\leq c\|u\|_{\dot{H}^{1.5}}^2\|q\|_{\dot{H}^1}.
$$
Using \eqref{eqaux} and $\alpha\geq0$, we obtain
$$
\frac{1}{2}\frac{d}{dt}\left(\|u\|_{\dot{H}^1}^2+\|\nabla\cdot{q}\|_{L^2}^2\right)+D_1\leq cD_1\|\nabla\cdot{q}\|_{L^2}^2.
$$
We conclude the result with the previous ideas.

\section*{Acknowledgment}
The author is funded by the Labex MILYON and the Grant MTM2014-59488-P from the Ministerio de Econom\'ia y Competitividad (MINECO, Spain).

\bibliographystyle{abbrv}

\end{document}